\documentclass[10pt]{amsart}
\usepackage{mathrsfs}
\usepackage{amsfonts} 
\usepackage{graphicx,latexsym,bm,amsmath,amssymb,verbatim,multicol,lscape}
\theoremstyle{definition}

\theoremstyle{remark}

\numberwithin{equation}{section}

\begin{document}
\title[New Lower Bounds for Least Common Multiples of Arithmetic Progressions]
{New Lower Bounds for the Least Common Multiples of Arithmetic Progressions}
\author{Rongjun Wu}
\address{Mathematical College, Sichuan University, Chengdu 610064, P.R. China}
\email{eugen\_woo@163.com}
\author{Qianrong Tan}
\address{School of Mathematics and Computer Science, Panzhihua University,
Panzhihua 617000, P.R. China}
\email{tqrmei6@126.com}
\author{Shaofang Hong*}
\address{Mathematical College, Sichuan University, Chengdu 610064, P.R. China and
Yangtze Center of Mathematics, Sichuan University, Chengdu 610064, P.R. China}
\email{sfhong@scu.edu.cn, s-f.hong@tom.com, hongsf02@yahoo.com }

\thanks{*Hong is the corresponding author and was supported partially by the National Science
Foundation of China Grant \# 10971145 and by the Ph.D. Programs Foundation of Ministry
of Education of China Grant \#20100181110073}

\keywords{Arithmetic progression; least common multiple; lower bound}
\subjclass[2000]{Primary 11B25, 11N13, 11A05}
\begin{abstract}
For relatively prime positive integers $u_0$ and $r$ and for $0\le k\le n$,
define $u_k:=u_0+kr$. Let $L_n:={\rm
lcm}(u_0, u_1, ..., u_n)$ and let $a, l\ge 2$ be any integers. In
this paper, we show that, for integers $\alpha \geq a$ and $r\geq
\max(a, l-1)$ and $n\geq l\alpha r$, we have $$L_n\geq
u_0r^{(l-1)\alpha +a-l}(r+1)^n.$$ Particularly, letting $l=3$
yields an improvement to the best previous lower bound on $L_n$
obtained by Hong and Kominers.
\end{abstract}

\maketitle

\section{Introduction}\label{intro}

Hanson and Nair initiated the search for effective estimates for the least
common multiple of the terms in a finite arithmetic progression; and, in
\cite{[Ha]} and in \cite{[N]} they managed to produce good upper
and lower bounds for ${\rm lcm}(1, 2, ..., n)$. In particular,
Nair \cite{[N]} discovered a nice new proof for the following well-known
nontrivial lower bound
\begin{equation}
{\rm lcm}(1, 2, ..., n)\ge 2^{n-1}
\end{equation}
for any integer $n\ge 1$. In \cite{[F3]}, Farhi provided an identity
involving the least common multiple of binomial coefficients and then use
it to give a simple proof of the estimate (1.1). Inspired
by Hanson's and Nair's works, Bateman, Kalb, and Stenger \cite{[BKS]}
and Farhi \cite{[F1]} respectively sought asymptotics and nontrivial lower
bounds for the least common multiples of arithmetic progressions. Recently,
Hong, Qian and Tan \cite{[HQT]} extended the Bateman-Kalb-Stenger theorem
from the linear polynomial to the product of linear polynomials. On the
other hand, Farhi \cite{[F1]} obtained several nontrivial bounds and posed a
conjecture which was later confirmed by Hong and Feng~\cite{[HF]}.
Hong and Feng~\cite{[HF]} also got an improved lower
bound for sufficiently long arithmetic progressions; this result
was later sharpened further by Hong and Yang \cite{[HY1]}. We notice that Hong
and Yang \cite{[HY2]} and Farhi and Kane \cite{[FK]} obtained
some related results regarding the least common multiple of a finite
number of consecutive integers. The theorem of Farhi and
Kane~\cite{[FK]} was extended by Hong and Qian~\cite{[HQ]}
from the set of positive integers to the general arithmetic progression case.
Recently, Qian, Tan and Hong \cite{[QTH]} obtained some results about the
least common multiple of consecutive terms in a quadratic progression.

In this paper, we study finite arithmetic progressions
$\{u_k:=u_0+kr\}_{k=0}^n$ with $u_0, r\geq 1$ being integers satisfying
$(u_0,r)=1$. Throughout, we define $L_n:={\rm lcm}(u_0,u_1,\ldots, u_n)$
to be the least common multiple of the sequence $\{u_k\}_{k=0}^n$.
We begin with the following lower bound on $L_n$:\\

{\bf Theorem 1.1.} \cite{[HY1]} \label{hythm} {\it Let $\alpha\geq 1$
be an integer. If $n>  r^\alpha$, then we have $L_n\geq u_0
r^\alpha(r+1)^n$.}\\

\noindent If $r=1$, then Theorem 1.1 is the
conjecture of Farhi~\cite{[F1]} proven by Hong and Feng~\cite{[HF]}.
If $\alpha =1$, then Theorem 1.1 becomes the improved lower
bound of Hong and Feng~\cite{[HF]}. In \cite{[HK]}, Hong and
Kominers sharpened the lower bound in Theorem 1.1 whenever
$\alpha , r\ge 2$. In particular, they proved the following theorem
which replaces the exponential condition $n>r^{\alpha }$ of
Theorem 1.1 with a linear condition $n\ge 2\alpha r$.\\

{\bf Theorem 1.2.} \label{mainthm} \cite{[HK]} {\it Let $a\ge 2$ be
any given integer. Then for any integers $\alpha , r\ge a$ and $n\ge
2\alpha r$, we have $L_n\ge u_0r^{\alpha +a-2}(r+1)^n$.}\\

Letting $a=2$, we see that Theorem 1.2 improves upon
Theorem 1.1 for all but three choices of $\alpha , r\ge 2$.
In the present paper, we provide a more general lower bound as
follows.\\

{\bf Theorem 1.3.} \label{mainthm1} {\it Let $a, l\ge 2$ be any
given integers. Then for any integers $\alpha \geq a$ and $r\geq
\max(a, l-1)$ and $n\geq l\alpha r$, we have
$L_n\ge u_0r^{(l-1)\alpha +a-l}(r+1)^n$.}\\

Picking $l=2$, then Theorem 1.3 becomes Theorem 1.2. Letting $l=3$
in Theorem 1.3 gives us the following new lower bound.\\

{\bf Theorem 1.4.} \label{mainthm} {\it Let $a\ge 2$ be any given
integer. Then for any integers $\alpha , r\ge a$ and $n\ge 3\alpha
r$, we have $L_n\ge u_0r^{2\alpha +a-3}(r+1)^n$.}\\

Since $\alpha \ge a\ge 2,$ we have $2\alpha +a-3>\alpha +a-2$.
Therefore the lower bound in Theorem 1.4 is better than that of
Theorem 1.2 when $n$ is large enough.

This paper is organized as follows. In Section \ref{sec2},
we first introduce relevant notation and previous results.
Finally, we prove Theorem 1.3.

\section{Proof of Theorem 1.3}\label{sec2}

For any real numbers~$x$ and~$y$, we say that~$y$
\emph{divides}\/~$x$ if there exists an integer~$z$ such that
$x=yz$. If $x$ divides $y$, then we write $y\mid x$. As usual,
we let $\lfloor x \rfloor$ denote the largest integer no more than
$x$.

Following Hong and Yang~\cite{[HY1]}, we denote, for each integer
$0\leq k\leq n$,
$$C_{n,k}:=\frac{u_k\cdots u_n}{(n-k)!},\quad L_{n,k}:={\rm lcm}(u_k,\ldots,u_n).$$
From the latter definition, we have that $L_n=L_{n,0}$.

The following Lemma first appeared in~\cite{[F1]} and was reproved
in \cite{[F2]} and \cite{[HF]}.\\

{\bf Lemma 2.1.} \cite{[F1]} \cite{[F2]} \cite{[HF]}\label{L1} {\it
For any integer $n\geq 1$, $C_{n,0}\mid L_n$.}\\

From Lemma 2.1, we see immediately that
\begin{equation}
\label{keq}L_{n,k}=A_{n,k}\frac{u_k\cdots u_n}{(n-k)!}=A_{n,k}\cdot
C_{n,k}
\end{equation}
for some integer $A_{n,k}\geq 1$.

Following Hong and Feng~\cite{[HF]} and Hong and Yang~\cite{[HY1]},
we define, for any $n\geq 1$,
\begin{equation}
k_n:=\max\left\{0,\left\lfloor
\frac{n-u_0}{r+1}\right\rfloor+1\right\}.\label{hat}
\end{equation}
Hong and Feng \cite{[HF]} proved the following result.\\

{\bf Lemma 2.2.} \cite{[HF]}\label{L2} {\it For all $n\geq 1$ and
$0\leq k\leq n$,}
$$
L_n\geq L_{n,k_n}\geq C_{n,k_n}\geq  u_0(r+1)^n.
$$

Now we are in a position to prove a lemma whose proof closely follows
the approach of Hong and Yang~\cite{[HY1]}.\\

{\bf Lemma 2.3.} \label{lemma2} {\it Let $a, l\ge 2$ be any given
integers. Then for any integers $\alpha \geq a$ and $r\geq \max(a,
l-1)$ and $n\geq l\alpha r$, we have $n-k_n>((l-1)\alpha +a-l)r$.}

\begin{proof}
If $n\le u_0$, then by the definition (2.2), $k_n\le 1$. Since
$\alpha , r\ge a\ge 2$ and $n\ge l\alpha r$, we derive that $n-k_n\ge
n-1\ge l\alpha r-1>((l-1)\alpha +a-l)r$.

Now we suppose that $n>u_0$. In this case, we have
$$k_n=\left\lfloor \frac{n-u_0}{r+1} \right\rfloor+1.$$
So we have
$$
k_n\le \frac{n-u_0}{r+1}+1\le \frac{n-1}{r+1}+1=\frac{n+r}{r+1}.
$$
It then follows that
\begin{equation}\label{part1}
n-k_n\ge n-\frac{n+r}{r+1}=\frac{(n-1)r}{r+1}\ge \frac{(l\alpha
r-1)r}{r+1}.
\end{equation}

Note that $r\ge l-1$ tells us that $r-l+1\ge 0$. Then from the
assumption $\alpha , r\ge a$ it follows that
\begin{align}
(l\alpha r-1)-(r+1)((l-1)\alpha +a-l)\nonumber&=(r-l+1)\alpha
-1-(r+1)(a-l)
\\&\ge a(r-l+1)-1-(r+1)(a-l)\label{part2}\\\nonumber&=l(r-a)+l-1>0.
\end{align}
Therefore by \eqref{part2}, we infer that
\begin{equation}\label{part3}
\frac{l\alpha r-1}{r+1}>(l-1)\alpha +a-l.
\end{equation}
The desired result then follows immediately from \eqref{part1} and
\eqref{part3}.
\end{proof}

Using the similar argument as that of Theorem 1.1, by Lemma 2.3 we can
now prove Theorem 1.3 as the conclusion of this paper.

\begin{proof}[Proof of Theorem 1.3.]

By hypothesis, we have $\alpha , r\ge a\ge 2$, $l\ge 2$ and $n\ge
l\alpha r$. It follows from Lemma 2.3 that $r^{(l-1)\alpha +a-l} \mid (n-k_n)!$.
Thus, we may express $(n-k_n)!$ in the form $r^{(l-1)\alpha +a-l} \cdot
B_n=(n-k_n)!$, with $B_n\geq 1$ being an integer. Letting $k=k_n$ in
\eqref{keq}, we find that
$$r^{(l-1)\alpha +a-l} \cdot B_n\cdot L_{n,k_n}=A_{n,k_n}\cdot u_{k_n}\cdots u_n.$$
It then follows that $r^{(l-1)\alpha +a-l} \mid A_{n,k_n}$, since
the requirement $(r,u_0)=1$ implies that $(r,u_k)=1$ for all $0\leq
k\leq n$.  Then, we get from \eqref{keq} and Lemma 2.2 that
$$L_{n,k_n}\geq r^{(l-1)\alpha +a-l} C_{n,k_n}\geq u_0r^{(l-1)\alpha +a-l}(r+1)^n.$$
Therefore the statement of Theorem 1.3 follows immediately.
The proof of Theorem 1.3 is complete.
\end{proof}

\end{document}